\documentclass[11pt]{article}

\usepackage[T2A]{fontenc}
\usepackage[cp1251]{inputenc}
\usepackage{amsfonts}
\usepackage{eufrak}
\usepackage{amssymb}
\usepackage{amsmath}



\begin{document}

\newtheorem{theorem}{Theorem}
\newtheorem{lemma}{Lemma}
\newtheorem{corollary}{Следствие}
\newtheorem{definition}{Опеределение}

\centerline{\textbf{On a characterisation theorem for
\text{\boldmath $a$}-adic solenoids}}

\vskip 1 cm

\centerline{\textbf{G.M. Feldman}}

\vskip 1 cm

\noindent B. Verkin Institute for Low Temperature Physics
 and Engineering of the National Academy \\ of Sciences of Ukraine,
47, Nauky ave, Kharkiv, 61103, Ukraine; e-mail:
feldman@ilt.kharkov.ua

\bigskip

 \makebox[20mm]{ }\parbox{125mm}{ \small According to the Heyde theorem the Gaussian distribution on the real line    is characterized by the symmetry
 of the conditional distribution of one linear form
 of   independent random variables given another.  We prove an analogue  of this theorem for  linear forms of two independent random variables taking values in an \text{\boldmath $a$}-adic solenoid $\Sigma_{{\text{\boldmath $a$}}}$ without elements of order 2, assuming that the characteristic functions of the random  variables do not vanish, and coefficients
 of the linear forms are topological automorphisms of
  $\Sigma_{{\text{\boldmath $a$}}}$.}

\bigskip

\bigskip


\bigskip

\centerline{\textbf{1.  Introduction}}

\bigskip

 By the well-known Heyde theorem the Gaussian distribution on the real line    is characterized by the symmetry of the conditional distribution of one linear form
 of   independent random variables given another (\cite{He}, see also \cite[\S\,13.4.1]{KaLiRa}). For two independent random variables this theorem can be formulated as follows.

\bigskip

\noindent{\bf Theorem A.}  {\it  Let $\xi_1$ and $\xi_2$ be independent random variables with distributions  $\mu_1$ and $\mu_2$.
Assume that
$\alpha\ne -1$. If the conditional distribution of the linear form
$L_2 = \xi_1 +
\alpha\xi_2$ given $L_1 = \xi_1 + \xi_2$ is symmetric,
then  $\mu_1$ and $\mu_2$ are Gaussian distributions.}

\bigskip

 Group analogs of Heyde's theorem in the case when independent random variables
take values in a locally compact Abelian group $X$,
and coefficients of the linear forms are topological automorphisms of
 $X$ were studied in the articles \cite{Fe2}--\cite{Fe3}, \cite{Fe20bb}--\cite{FeTVP}, \cite{My2}--\cite{MiFe1}, see also   \cite[Chapter VI]{Fe5}). We remark that in all cited
 articles the corresponding characterization theorems were proved under certain restrictions
 on
 coefficients of the linear forms. In the article we prove an analogue of Heyde's theorem for
 two independent random variables taking values in an  \text{\boldmath $a$}-adic solenoid
 without elements of order 2, assuming that the characteristic function of considering random variables do not vanish.
 It is important that we do not impose any restrictions on coefficients of the linear forms.

 We note that although our proofs of the main results use  methods of abstract harmonic analysis, it was quite unexpected that some
facts of complex analysis were also used, in particular, the well known Hadamard theorem on the representation of an entire function of finite order.

Before we formulate the main theorem recall some definitions
and agree on notation.  Let $X$ be a second countable locally compact Abelian group. We will consider only such groups, without mentioning it specifically. Denote by ${\rm Aut}(X)$ the group
of topological automorphisms of $X$, and by  $I$ the identity automorphism of a group. Denote by $Y$ the character
group of the group $X$, and by  $(x,y)$ the value of a character $y \in Y$ at an element $x \in X$. If $G$ is a closed subgroup of $X$, denote by
 $A(Y, G) = \{y \in Y: (x, y) = 1$ \mbox{ for all } $x \in G \}$
its annihilator.    Let $X_1$ and
$X_2$ be locally compact Abelian groups with character groups
$Y_1$ and $Y_2$ respectively. Let
 $\alpha:X_1\mapsto X_2$
be a continuous homomorphism.
The adjoint homomorphism $\tilde\alpha: Y_2\mapsto Y_1$
is defined by the formula $(\alpha x_1,
y_2)=(x_1 , \tilde\alpha y_2)$ for all $x_1\in X_1$, $y_2\in
Y_2$.   Denote by $\mathbb{C}$   the set of complex numbers, by $\mathbb{R}$ the group of real numbers, by
 $\mathbb{Z}$ the group of integers,  by ${\mathbb Z}(m)=\{0, 1, \dots, m-1\}$ the group of  residue  classes modulo $m$, and by $\mathbb{T}=\{z\in \mathbb{C}: |z|=1\}$ the circle group.
 Let $f(y)$ be a function on the group    $Y$,   and let $h \in
Y$. Denote by   $\Delta_h$   the finite difference operator
$$
\Delta_h f(y) = f(y + h) - f(y), \quad y \in Y.
$$
Let  $n$ be an integer. Denote by $f_n:X \mapsto X$ an  endomorphism of the group $X$
 defined by the formula  $f_nx=nx$, $x\in X$.
Put $X^{(n)}=f_n(X)$, $X_{(n)}={\rm Ker}f_n.$

Let $\mu$ be a measure or a signed measure on the group $X$. The characteristic function (the  Fourier transform)
of   $\mu$ is defined by the formula
$$
\hat\mu(y) =
\int_{X}(x, y)d \mu(x), \quad y\in Y.
$$
Denote by ${\rm M}^1(X)$ the
convolution semigroup of probability measures (distributions) on the group $X$. Let
${\mu\in {\rm M}^1(X)}$.
Denote by $\sigma(\mu)$ the support of $\mu$. Define
the distribution $\bar \mu \in {\rm M}^1(X)$ by the formula
 $\bar \mu(B) = \mu(-B)$ for any Borel subset $B$ of $X$.
Note that $\hat{\bar{\mu}}(y)=\overline{\hat\mu(y)}$. If $G$ is a Borel subgroup of $X$, we
denote by ${\rm M}^1(G)$ the subsemigroup of  ${\rm M}^1(X)$  of distributions concentrated
on $G$.
A distribution  $\gamma\in {\rm M}^1(X)$  is called Gaussian
(see \cite[Chapter IV, \S 6]{Pa})
if its characteristic function is represented in the form
\begin{equation}\label{f1}
\hat\gamma(y)= (x,y)\exp\{-\varphi(y)\}, \quad  y\in Y,
\end{equation}
where $x \in X$, and $\varphi(y)$ is a continuous non-negative function
on the group $Y$
 satisfying the equation
 \begin{equation}\label{f2}
\varphi(u + v) + \varphi(u
- v) = 2[\varphi(u) + \varphi(v)], \quad u,  v \in
Y.
\end{equation}
 Note that, in particular, the degenerate distributions are Gaussian.
Denote by $\Gamma(X)$ the set of Gaussian distributions on
    $X$. Denote by   $E_x$  the degenerate distribution
 concentrated at a point $x\in X$.

\bigskip

\centerline{\textbf{2.   \text{\boldmath $a$}-adic solenoids}}
\bigskip

 Recall the definition of an $\text{\boldmath $a$}$-adic solenoid.
Put \text{\boldmath $a$}= $(a_0, a_1,\dots)$, where all $a_j \in {\mathbb{Z}}$,
$a_j > 1$. Denote by $\Delta_{{\text{\boldmath $a$}}}$
 the group of  ${\text{\boldmath $a$}}$-adic integers.  Consider the group
$\mathbb{R}\times\Delta_{{\text{\boldmath $a$}}}$. Denote by
 $B$ the subgroup of the group
 $\mathbb{R}\times\Delta_{{\text{\boldmath $a$}}}$ of the form
$B=\{(n,n\mathbf{u})\}_{n=-\infty}^{\infty}$, where
$\mathbf{u}=(1, 0,\dots,0,\dots)$. The factor-group $\Sigma_{{\text{\boldmath $a$}}}=(\mathbb{R}\times\Delta_{{\text{\boldmath $a$}}})/B$ is called
  an ${\text{\boldmath $a$}}$-{adic solenoid}.   The group $\Sigma_{{\text{\boldmath $a$}}}$ is  compact, connected and has
dimension 1   (\cite[(10.12), (10.13),
(24.28)]{Hewitt-Ross}).
 The character group of the group   $\Sigma_{{\text{\boldmath $a$}}}$ is topologically isomorphic
to the discrete group of the form
$$H_{\text{\boldmath $a$}}= \left\{{m \over a_0a_1 \dots a_n} : \ n = 0, 1,\dots; \ m
\in {\mathbb{Z}} \right\}.
$$
 In order not to complicate the notation we will identify $H_{\text{\boldmath $a$}}$ with the character group of the group $\Sigma_{{\text{\boldmath $a$}}}$. We will also consider  $H_{\text{\boldmath $a$}}$ as a subset of  $\mathbb{R}$. Any topological automorphism $\alpha$
 of the group
$\Sigma_\text{\boldmath $a$}$  is of the form
$\alpha = f_p f_q^{-1}$ for some mutually prime $p$ and $q$,
where $f_p, f_q \in {\rm
Aut}(\Sigma_\text{\boldmath $a$})$.
  We will identify $\alpha = f_p f_q^{-1}$
with the real number ${p\over q}$. If   $\alpha = f_p f_q^{-1}$, then $\tilde\alpha=f_p f_q^{-1}$, and we will also identify $\tilde\alpha$ with the real number ${p\over q}$.
We note that if $\alpha\ne-I$, then either ${\rm Ker}(I+\alpha)=\{0\}$  or ${\rm Ker}(I+\alpha)\cong \mathbb{Z}(m)$ for some $m$.
Put $G=(\Sigma_\text{\boldmath $a$})_{(2)}$. It is easy to verify that there are only two possibilities for $G$:
 either $G=\{0\}$  or
  $G\cong \mathbb{Z}(2)$.
  It is obvious that $G=\{0\}$ if and only if  $f_2\in {\rm Aut}(\Sigma_\text{\boldmath $a$})$.

  It follows from  (\ref{f1}) и (\ref{f2}) that the characteristic function of a Gaussian distribution   $\gamma$ on  an   \text{\boldmath $a$}-adic solenoid
 $\Sigma_\text{\boldmath $a$}$ is of the form
$$
\hat\gamma(y)=(x, y)\exp\{-\sigma y^2\}, \quad y\in H_{\text{\boldmath $a$}},
$$
where $x\in  \Sigma_\text{\boldmath $a$}$, $\sigma \ge 0$.

Let $\xi_1$ and $\xi_2$ be independent random variables with values in a locally compact Abelian  group
       $X$  and distributions
  $\mu_1$ and $\mu_2$. Let $\alpha_j, \beta_j\in{\rm Aut}(X)$. Consider the linear forms $L_1=\alpha_1\xi_1+\alpha_2\xi_2$ and $L_2=\beta_1\xi_1+\beta_2\xi_2$ and assume that
the conditional distribution of the linear form $L_2$ given $L_1$ is symmetric.
It is easy to see that the description of possible distributions $\mu_j$ is reduced to the case when  $L_1=\xi_1+\xi_2$, $L_2=\xi_1+\alpha\xi_2$, where $\alpha\in{\rm Aut}(X)$.

The main result of the article is the following theorem.

\begin{theorem}.\label{th1}  {\it  Consider an \text{\boldmath $a$}-adic solenoid $X=\Sigma_\text{\boldmath $a$}$. Assume that $X$ contains no elements
of order $2$.
Let $\alpha$ be a topological automorphism of the group $X$.
 Put $K={\rm Ker}(I+\alpha)$.
Let $\xi_1$ and $\xi_2$ be independent random variables with values in the group
       $X$  and distributions
  $\mu_1$ and $\mu_2$ with nonvanishing characteristic functions. Assume that the conditional distribution of the linear form
$L_2=\xi_1+\alpha\xi_2$
 given $L_1 = \xi_1 + \xi_2$ is symmetric.
Then each of the distributions $\mu_j$ can be represented in the form
$\mu_j=\gamma_j*\omega$, where $\gamma_j\in \Gamma(X)$, $\omega\in {\rm M^1}(K)$. Moreover, if $\alpha>0$, then $\mu_j=E_{x_j}*\omega$, where
$x_j\in X$, $j=1, 2$.}
\end{theorem}

\noindent{\bf Corollary 1.} {\it  Consider an \text{\boldmath $a$}-adic solenoid $X=\Sigma_\text{\boldmath $a$}$. Assume that $X$ contains no elements
of order $2$.
Let $\alpha$ be a topological automorphism of the group $X$. Let $\xi_1$ and $\xi_2$ be
independent random variables with values in the group
       $X$  and distributions
  $\mu_1$ and $\mu_2$ with nonvanishing characteristic functions. The symmetry of the
  conditional distribution of the linear form   $L_2 = \xi_1 + \alpha\xi_2$ given $L_1 = \xi_1
  +
\xi_2$ implies that  $\mu_j$ are Gaussian distributions if and only if  ${\rm Ker}(I+\alpha)=\{0\}$.}

\bigskip

 Theorem 1 will be proved in \S 3. In \S 4 we prove that generally speaking, Theorem 1 fails
 if
 an \text{\boldmath $a$}-adic solenoid $X=\Sigma_\text{\boldmath $a$}$ contains an element of
 order 2.
Namely, we will prove for such groups that if   $\alpha<0$, then
there exist independent
 random variables $\xi_1$ and $\xi_2$  with values in the group
       $X$ and distributions  $\mu_j\notin\Gamma(X)*{\rm M}^1(K)$ with nonvanishing
       characteristic functions such that   the conditional distribution of the linear
       form
$L_2=\xi_1+\alpha\xi_2$
 given $L_1 = \xi_1 + \xi_2$ is symmetric.

\bigskip

\centerline{\textbf{3.  Proof of Theorem 1}}

\bigskip

 To prove Theorem 1 we need some lemmas.

\begin{lemma}\label{lem1}{\rm (\cite[Lemma 16.1]{Fe5})}. {\it
Let $X$ be a locally compact
Abelian group, $Y$ be
its character group, $\alpha$ be a topological automorphism of $X$.
Let $\xi_1$ and $\xi_2$ be independent random variables with values in the group
       $X$  and distributions
 $\mu_1$ and $\mu_2$. The conditional distribution of the linear form
$L_2 = \xi_1 + \alpha\xi_2$
 given $L_1 = \xi_1 + \xi_2$ is symmetric if and only
 if the characteristic functions $\hat\mu_j(y)$  satisfy Heyde's functional equation}
\begin{equation}\label{42}
\hat\mu_1(u+v )\hat\mu_2(u+\tilde\alpha v )=
\hat\mu_1(u-v )\hat\mu_2(u-\tilde\alpha v), \quad u, v \in Y.
\end{equation}
\end{lemma}

It is convenient for us to formulate as   lemmas the following well-known
statements  (see e.g.   \cite[\S 2]{Fe5}).
\begin{lemma}\label{lem2}. {\it Let $X$ be a locally compact
Abelian group,   $Y$ be its character group. Let $\mu\in{\rm
M}^1(X)$. Then the set $E=\{y\in Y:\ \hat\mu(y)=1\}$ is a closed subgroup
of $Y$, the characteristic function
$\hat\mu(y)$ is $E$-invariant, i.e. $\hat\mu(y)$ takes a constant value on each coset
of $E$ in the
group $Y$, and $\sigma(\mu)\subset A(X,E)$.}
\end{lemma}
\begin{lemma}\label{lem3}. {\it Let $X$ be a locally compact Abelian group,
 $G$ be a Borel subgroup of  $X$,
$\mu\in {\rm M}^1(G)$, $\mu=\mu_1*\mu_2$, where $\mu_j\in {\rm
M}^1(X)$. Then the distributions $\mu_j$ can be replaced by their
shifts $\mu'_1=\mu_1*E_x$, $\mu'_2=\mu_2*E_{-x}$, $x\in X$, in such a manner that
$\mu'_j\in {\rm M}^1(G)$, $j=1, 2$.}
\end{lemma}

 We will formulate as a lemma the following easily verifiable statement   (see e.g.
 \cite[Lemma 6.9]{Fe9}).

\begin{lemma}\label{lem7}.   {\it Let $X=\mathbb{R}\times G$, where $G$ is a locally
compact Abelian group, $Y$ and $H$ be  the character groups of the groups
$X$  and $G$ respectively. Denote by $(t, g)$, $t\in\mathbb{R}$, $g\in G$, elements of
the group $X$ and by  $(s, h)$, $s\in\mathbb{R}$, $h\in H$, elements of the group  $Y$. Let
$\mu\in {\rm M}^1(X)$ and assume that $\hat\mu(s, 0)$, $s\in \mathbb{C}$, is
an entire
function in $s$. Then $\hat\mu(s, h)$ is an entire function in
$s$ for every fixed   $h\in H$,
the representation
$$
\hat\mu(s, h)=\int_X\exp\{its\}(g, h)d\mu(t, g)
$$
holds for all $s\in \mathbb{C}$, $h\in H$, and the inequality
\begin{equation}\label{8_16}
\max_{|s|\le r}|\hat\mu(s, h)|\le \max_{|s|\le r}|\hat\mu(s, 0)|, \quad h\in H,
\end{equation}
is valid. Moreover, the function
${\hat\mu(-iy+x, h)/\hat\mu(-iy, 0)}$, where $x, y\in \mathbb{R}$, for any fixed $y$ is
a characteristic function of variable $(x, h)\in \mathbb{R}\times H$.}
\end{lemma}

\noindent{\bf Proof of theorem  1.} Denote by  $Y$ the character group of the group $X$, and consider
$Y$ as a subset of $\mathbb{R}$.
Let $\alpha = f_p f_q^{-1}$ for some mutually prime  $p$ and $q$,
where $f_p, f_q \in {\rm
Aut}(X)$.
By Lemma \ref{lem1},  the symmetry of the conditional distribution of the linear form
$L_2$
 given $L_1 $ implies that the characteristic functions $\hat\mu_j(y)$  satisfy  Heyde's
 functional equation $(\ref{42})$. Put $\nu_j = \mu_j* \bar \mu_j$.
Then  $\hat \nu_j(y) = |\hat \mu_j(y)|^2 > 0,$  $y \in Y$.
Obviously, the characteristic functions $\hat \nu_j(y)$
also satisfy Heyde's
 functional equation $(\ref{42})$, which takes the form
\begin{equation}\label{20_04_3}
\hat\nu_1(u+v )\hat\nu_2(u+\tilde\alpha v )=
\hat\nu_1(u-v )\hat\nu_2(u-\tilde\alpha v), \quad u, v \in Y.
\end{equation}

Describe the scheme of the proof of Theorem 1. First we get some representation
 for the characteristic functions $\hat \nu_j(y)$.
By the topological automorphism $\alpha$ we find natural numbers $m$ and $n$ and consider the group $\mathbb{R}\times \mathbb{Z}(mn)$. Next we
  construct distributions $M_j\in{\rm M}^1(\mathbb{R}\times \mathbb{Z}(mn))$ and a continuous monomorphism
${\pi:\mathbb{R}\times \mathbb{Z}(mn)\mapsto X}$ such that
 $\nu_j=\pi(M_j)$. In so doing the characteristic functions of the distributions
  $M_j$ satisfy some Heyde's functional equation.   Next we solve  the obtained Heyde
  functional
 equation  and receive the representation for the characteristic functions of the distributions $M_j$. Finally, we find
 the desired
representation for the  distributions $\mu_j$.

First prove the theorem assuming that  $\alpha\ne \pm I$.
Put $\varphi_j(y) =  \log \hat\nu_j(y)$, $j=1, 2$.
It follows from $(\ref{20_04_3})$ that the functions $\varphi_j(y)$
satisfy the equation
\begin{equation}\label{12}
\varphi_1(u+v) + \varphi_2(u+\tilde\alpha v)- \varphi_1(u-v)
- \varphi_2(u-\tilde\alpha v)=0, \quad u, v \in Y.
  \end{equation}
We use the finite difference method to solve this equation.
Let $k_1$  be an arbitrary element of  $Y$. Put $h_1=\tilde\alpha k_1$. Hence,
$h_1 -\tilde\alpha k_1 = 0$.
Substitute in
  (\ref{12}) $u+h_1$ for $u$ and $v+k_1$ for $v$. Subtracting equation
  (\ref{12}) from the resulting equation, we obtain
\begin{equation}\label{13}
\Delta_{l_{11}}\varphi_1(u+v) + \Delta_{l_{12}}\varphi_2(u+\tilde\alpha v) -\Delta_{l_{13}}\varphi_1(u-v)=0, \quad u, v \in
Y,
\end{equation}
where $l_{11}= (I+\tilde\alpha)k_1$, $l_{12}=2 \tilde\alpha k_1$, $l_{13}= (\tilde\alpha-I)k_1.$

 Let  $k_2$  be an arbitrary element of the group $Y$. Put $h_2=k_2$.
Hence, $h_2-k_2 = 0$. Substitute
in   (\ref{13}) $u+h_2$  for $u$ and $v+k_2$ for $v$.
Subtracting equation
  (\ref{13}) from the resulting equation, we find
\begin{equation}\label{14}
\Delta_{l_{21}}\Delta_{l_{11}}\varphi_1(u+v) +
\Delta_{l_{22}}\Delta_{l_{12}}\varphi_2(u+\tilde\alpha v) =
0,  \quad u, v \in Y,
\end{equation}
where $l_{21}=2 k_2$, $l_{22}=(I+\tilde\alpha)k_2.$

Let $k_3$ be an arbitrary element of $Y$. Put $h_3 =-\tilde\alpha
k_3$. Hence, $h_3 +\tilde\alpha k_3 = 0$. Substitute in
(\ref{14}) $u+h_3$ for $u$ and $v+k_3$ for $v$.   Subtracting equation
  (\ref{14}) from the resulting equation, we get
\begin{equation}\label{15}
\Delta_{l_{31}}\Delta_{l_{21}}\Delta_{l_{11}}\varphi_1(u+v)
= 0, \quad u, v \in Y,
\end{equation}
where $l_{31}= (I-\tilde\alpha)k_3.$ Putting
in (\ref{15}) $v=0$, we find
\begin{equation}\label{16}
\Delta_{l_{31}}\Delta_{l_{21}}\Delta_{l_{11}}\varphi_1(u) = 0,
 \quad u \in Y.
\end{equation}
 Reasoning similarly we obtain from  (\ref{14}) that the function $\varphi_2(y)$
satisfies the equation
\begin{equation}\label{new16}
\Delta_{l_{32}}\Delta_{l_{22}}\Delta_{l_{12}}\varphi_2(u) = 0,
 \quad u \in Y,
\end{equation}
where $l_{32}=-(I-\tilde\alpha)k_3$.
Put \begin{equation}\label{20_08_1}
H=(I+\tilde\alpha)(Y)\cap Y^{(2)}\cap(I-\tilde\alpha)(Y).
\end{equation}
It follows from (\ref{16}) and (\ref{new16}) that the functions
 $\varphi_j(y)$ satisfy the equation
\begin{equation}\label{17}
\Delta_h^3\varphi_j(y)=0, \quad  h\in H, \quad y\in Y, \quad j=1, 2.
\end{equation}
It is useful to remark that we received equation (\ref{17})
without using the fact that the group $X$
contains no elements of order 2.

Since the group $X$ contains no elements of order $2$, we have $Y^{(2)}=Y$.  It is obvious
that
$(I+\tilde\alpha)(Y)=Y^{(m)}$ for some natural $m$. We can assume
that $m$ is is the minimum possible. Similarly,
$(I-\tilde\alpha)(Y)=Y^{(n)}$ for some natural $n$, and we can   also assume
that $n$ is  the minimum possible.
It follows from $Y^{(2)}=Y$ that $m$ and $n$ are odd.
Since $m$ and $n$ are divisors of $p+q$ and
$p-q$ respectively and $p$ and $q$ mutually prime,  $m$ and $n$ are also
mutually prime.
 This implies that $H=Y^{(m)}\cap Y^{(n)}=Y^{(mn)}$.   Consider the factor-group $Y/H$. It is obvious that
 $Y/H\cong \mathbb{Z}(mn)$.
 Take an element
 $y_0\in Y$ in such a way that the coset $y_0+H$ be  a generator of the factor-group $Y/H$. Then
\begin{equation}\label{23_04_1}
Y=H\cup(y_0+H)\cup(2y_0+H)\cup\dots\cup((mn-1)y_0+H)
\end{equation}
is a decomposition of the group $Y$ with respect to the subgroup $H$.
Put
$$\psi_{lj}(y)=\varphi_j(ly_0+y), \quad y\in H, \quad l=0, 1,\dots, mn-1, \quad j=1, 2.
$$
It follows from  (\ref{17}) that the functions $\psi_{lj}(y)$ satisfy the equation
$$
\Delta_h^3\psi_{lj}(y)=0, \quad h, y\in H, \quad l=0, 1,\dots, mn-1, \quad j=1, 2.
$$
This implies that for any coset $ly_0+H$ there exists the polynomial
\begin{equation}\label{13_05_1}
P_{lj}(y)=A_{lj}y^2+B_{lj}y+C_{lj}, \quad j=1, 2,
\end{equation}
on $\mathbb{R}$ with real coefficients such that
$$
P_{lj}(y)=\varphi_j(y), \quad y\in ly_0+H, \quad l=0, 1,\dots, mn-1, \quad j=1, 2.
$$
It follows from this that
\begin{equation}\label{20_04_1}
\hat\nu_{j}(y)=\exp\{A_{lj}y^2+B_{lj}y+C_{lj}\}, \quad      y\in ly_0+H, \quad l=0, 1,\dots, mn-1, \quad j=1, 2.
\end{equation}

Consider the group $\mathbb{R}\times \mathbb{Z}(mn)$ and denote by $(s, l)$, $s\in \mathbb{R}$, $l\in \mathbb{Z}(mn)$ its elements. Define the mapping $\tau:Y\mapsto \mathbb{R}\times \mathbb{Z}(mn)$ by the following way:
$\tau y=(y,  l)$  if $y\in ly_0+H$, $l=0, 1,\dots, mn-1$.  Obviously,
$\tau$ is a monomorphism and the subgroup $\tau(Y)$ is
dense in $\mathbb{R}\times \mathbb{Z}(mn)$. Since $m$ and $n$ are divisors of
$p+q$ and
$p-q$ respectively, and $p$ and $q$ are mutually prime,     $p$ is mutually prime with $m$ and
with $n$, and $q$ is mutually prime with
  $m$ and with $n$. Hence, $f_p, f_q\in {\rm Aut}(\mathbb{Z}(mn))$.

Let $y\in ly_0+H$. If we consider $l$ as an element of the group  $\mathbb{Z}(mn)$,
and $f_p$ and $f_q$ as automorphisms of the group
 $\mathbb{Z}(mn)$, then it is easy to see that
$\tilde\alpha y\in (f_pf_q^{-1}l)y_0+H$. Therefore, on the one hand, $\tau\tilde\alpha y=(\tilde\alpha y, f_pf_q^{-1}l)=(f_pf_q^{-1} y, f_pf_q^{-1}l)$. On the other hand, since $\tau y=(y, l)$ and $f_pf_q^{-1}(y, l)=(f_pf_q^{-1} y, f_pf_q^{-1}l)$, we have
\begin{equation}\label{20_04_5}
\tau\tilde\alpha y=f_pf_q^{-1}\tau y, \quad y\in Y.
\end{equation}

Define on the set $\tau(Y)$ the functions $g_j(s, l)$ in the following way.
Let $(s, l)=\tau y\in \tau(Y)$. Put $g_j(s, l)=\hat\nu_j(y)$, $j=1, 2$.
Since $\tau$ is a monomorphism, the functions $g_j(s, l)$ are correctly defined.
It follows from (\ref{20_04_1}) that the functions $g_j(s, l)$ are   continuous on the subgroup $\tau(Y)$ in the topology induced on  $\tau(Y)$ by the topology of the group
$\mathbb{R}\times \mathbb{Z}(mn)$, and taking into account that the subgroup $\tau(Y)$ is dense
in $\mathbb{R}\times \mathbb{Z}(mn)$, the functions  $g_j(s, l)$ can be extended by continuity
to some
continuous positive definite functions  $\bar g_j(s, l)$ on the group $\mathbb{R}\times \mathbb{Z}(mn)$. By the Bochner theorem,  there exist the distributions $M_j\in{\rm M}^1(\mathbb{R}\times \mathbb{Z}(mn))$  such that
$\hat M_j(s, l)=\bar g_j(s, l)$, $(s, l)\in \mathbb{R}\times \mathbb{Z}(mn)$. It is obvious that
\begin{equation}\label{20_04_2}
\hat M_j(\tau y)=\hat\nu(y), \quad y\in Y, \quad  j=1, 2.
\end{equation}
It follows from (\ref{20_04_1}) and (\ref{20_04_2}) that the characteristic functions
$\hat M_j(s, l)$  do not vanish, and hence
 $\hat M_j(s, l)>0$, $j=1, 2$.

Note that the character group of the group  $\mathbb{R}\times \mathbb{Z}(mn)$ is topologically isomorphic to $\mathbb{R}\times \mathbb{Z}(mn)$. In order not to complicate the notation we will assume that it coincides with
 $\mathbb{R}\times \mathbb{Z}(mn)$. Put $\pi=\tilde \tau$, $\pi:\mathbb{R}\times \mathbb{Z}(mn)\mapsto X$.
Since the subgroup $\tau(Y)$ is dense in $\mathbb{R}\times \mathbb{Z}(mn)$,  $\pi$ is a
continuous monomorphism generating an isomorphism of the semigroups of distributions  ${\rm M^1}(\mathbb{R}\times \mathbb{Z}(mn))$ and ${\rm M^1}(\pi(\mathbb{R}\times \mathbb{Z}(mn))).$
We shall also denote this isomorphism by $\pi$.
 It follows from (\ref{20_04_2}) that
 \begin{equation}\label{19_05_1}
\nu_j=\pi(M_j), \quad j=1, 2.
\end{equation}

Find a representation for the characteristic functions
  $\hat M_j(s, l)$.
  It follows from  (\ref{20_04_3}) and (\ref{20_04_2}) that
\begin{equation}\label{20_04_21}
\hat M_1(\tau u+\tau v )\hat M_2(\tau u+\tau\tilde\alpha v )=
\hat M_1(\tau u-\tau v )\hat M_2(\tau u-\tau\tilde\alpha v ), \quad u, v \in Y.
\end{equation}
Taking into account (\ref{20_04_5}), we find from (\ref{20_04_21}) that
\begin{equation}\label{20_04_6}
\hat M_1(\tau u+\tau v )\hat M_2(\tau u+f_pf_q^{-1}\tau v )=
\hat M_1(\tau u-\tau v )\hat M_2(\tau u-f_pf_q^{-1}\tau v ), \quad u, v \in Y.
\end{equation}
 Since the subgroup $\tau(Y)$ is dense in $\mathbb{R}\times \mathbb{Z}(mn)$,
 (\ref{20_04_6}) implies that the characteristic functions  $\hat M_j(s, l)$
on the group $\mathbb{R}\times \mathbb{Z}(mn)$ satisfy the following Heyde's
 functional equation
\begin{equation}\label{20_04_7}
\hat M_1(s_1+s_2, l_1+l_2 )\hat M_2(s_1+ pq^{-1}s_2, l_1+ f_pf_q^{-1}l_2 )$$$$
=\hat M_1(s_1-s_2, l_1-l_2 )\hat M_2(s_1-pq^{-1}s_2, l_1- f_pf_q^{-1}l_2 ), \quad (s_j, l_j)\in \mathbb{R}\times \mathbb{Z}(mn).
\end{equation}
Note that $\mathbb{Z}(mn)=\mathbb{Z}(m)\times \mathbb{Z}(n)$, and represent an element
$l\in \mathbb{Z}(mn)$ in the form $l=(a, b)$, $a\in \mathbb{Z}(m)$, $b\in \mathbb{Z}(n)$.
Since $m$ is a divisor of $p+q$, and  $n$ is a divisor of
$p-q$, it is easy to see that the automorphism   $f_pf_q^{-1}$ acts
on the group
$\mathbb{Z}(m)\times \mathbb{Z}(n)$ in the following way
$f_pf_q^{-1}(a, b)=(-a, b)$, and we can write equation   (\ref{20_04_7}) in the form
\begin{equation}\label{21_04_1}
\hat M_1(s_1+s_2, a_1+a_2, b_1+b_2 )\hat M_2(s_1+pq^{-1}s_2, a_1- a_2, b_1+b_2 )$$$$
=\hat M_1(s_1-s_2, a_1-a_2, b_1-b_2 )\hat M_2(s_1-pq^{-1}s_2, a_1+a_2, b_1-b_2 ),
\quad (s_j, a_j, b_j)\in \mathbb{R}\times \mathbb{Z}(m)\times \mathbb{Z}(n).
\end{equation}
Putting in (\ref{21_04_1})  $s_1=s_2=0$, $a_1=a_2=0$, we get
\begin{equation}\label{21_04_2}
\hat M_1(0, 0, b_1+b_2 )\hat M_2(0, 0, b_1+b_2 )
=\hat M_1(0, 0, b_1-b_2 )\hat M_2(0,  0, b_1-b_2 ), \quad b_j\in\mathbb{Z}(n).
\end{equation}
Putting in (\ref{21_04_2}) $b_1=b_2=b$, we get $\hat M_j(0, 0, 2b)=1$ for any
$b\in\mathbb{Z}(n)$, and taking into account that
 $n$ is odd, we conclude that $\hat M_j(0, 0, b)=1$, $j=1, 2,$
for $b\in\mathbb{Z}(n)$. By Lemma \ref{lem2}, this implies that
$\hat M_j(s, a, b)=\hat M_j(s, a, 0)$ for all
$(s, a, b)\in \mathbb{R}\times \mathbb{Z}(m)\times \mathbb{Z}(n)$  and $\sigma(M_j)\subset A(\mathbb{R}\times \mathbb{Z}(m)\times \mathbb{Z}(n), \mathbb{Z}(n))=\mathbb{R}\times \mathbb{Z}(m)$, $j=1, 2.$ Therefore we can write the characteristic function $\hat M_j(s, a, b)$
in the form $\hat M_j(s, a)$. Then equation  (\ref{21_04_1}) takes the form
\begin{equation}\label{21_04_3}
\hat M_1(s_1+s_2, a_1+a_2)\hat M_2(s_1+pq^{-1}s_2, a_1- a_2)$$$$
=\hat M_1(s_1-s_2, a_1-a_2)\hat M_2(s_1-pq^{-1}s_2, a_1+a_2),
\quad (s_j, a_j)\in \mathbb{R}\times \mathbb{Z}(m).
\end{equation}
Put in (\ref{21_04_3})  $a_1=a_2=0$. Taking into account Lemma \ref{lem1}, we get by Theorem
A from the obtaining equation that
\begin{equation}\label{11_05_3}
\hat M_j(s, 0)=\exp\{-\sigma_j s^2\}, \quad  s \in \mathbb{R}, \quad j=1, 2.
\end{equation}
By Lemma \ref{lem7},  $\hat M_j(s, a)$ are entire functions for any $a\in \mathbb{Z}(m)$.
As noted above the functions $\hat M_j(s, a)$ do not vanish. It follows from (\ref{8_16}) and (\ref{11_05_3}) that the entire function
$\hat M_j(s, a)$ is of at most  order 2.
Applying the Hadamard theorem on the representation of an entire function of finite order, we obtain that
\begin{equation}\label{12_05_1}
\hat M_j(s, a)=\exp\{\alpha_{aj}s^2+\beta_{aj}s+\gamma_{aj}\}, \quad  s \in \mathbb{R},  \ a\in \mathbb{Z}(m), \quad j=1, 2.
\end{equation}
Put
\begin{equation}\label{12_05_2}
\phi_j(s, a) = \alpha_{aj}s^2+\beta_{aj}s+\gamma_{aj},   \quad  s \in \mathbb{R},  \ a\in \mathbb{Z}(m), \quad j=1, 2.
\end{equation}
It follows from  $(\ref{21_04_3})$ that the functions $\phi_j(s, a)$
satisfy the eqaution
\begin{equation}\label{21_04_4}
\phi_1(s_1+s_2, a_1+a_2) + \phi_2(s_1+pq^{-1}s_2, a_1- a_2)$$$$- \phi_1(s_1-s_2, a_1-a_2)
-\phi_2(s_1-pq^{-1}s_2, a_1+a_2)=0, \quad (s_j, a_j)\in \mathbb{R}\times \mathbb{Z}(m).
  \end{equation}
Substituting  (\ref{12_05_2}) into (\ref{21_04_4}) and putting in the obtained equation $a_1=a_2=a$, we get
\begin{equation}\label{22_04_2}
\alpha_{2a1}(s_1+s_2)^2+\beta_{2a1}(s_1+s_2)+\gamma_{2a1}+
\alpha_{02}(s_1+pq^{-1}s_2)^2+\beta_{02}(s_1+pq^{-1}s_2)+\gamma_{02}$$$$-
\alpha_{01}(s_1-s_2)^2-\beta_{01}(s_1-s_2)-\gamma_{01}-
\alpha_{2a2}(s_1-pq^{-1}s_2)^2-\beta_{2a2}(s_1-pq^{-1}s_2)-\gamma_{2a2}=0.
\end{equation}
Equating  the coefficients  of
$s_1^2$ and  $s_2^2$ on each  side  of (\ref{22_04_2}), we get
\begin{equation}\label{18_05_1}
\alpha_{2a1}+\alpha_{02}-\alpha_{01}-\alpha_{2a2}=0, \quad \alpha_{2a1}+(pq^{-1})^2\alpha_{02}-\alpha_{01}-(pq^{-1})^2\alpha_{2a2}=0.
\end{equation}
Taking into account that $(pq^{-1})^2\ne 1$ and the fact that $m$ is odd, (\ref{18_05_1}) and (\ref{11_05_3})
imply that
 $\alpha_{a1}=\alpha_{01}=\sigma_1$ and $\alpha_{a2}=\alpha_{02}=\sigma_2$ for all $a\in \mathbb{Z}(m)$.
Equating  the coefficients  of $s_1$ and $s_2$ on each side  of (\ref{22_04_2}) we obtain
\begin{equation}\label{22_04_3}
\beta_{2a1}+\beta_{02}-\beta_{01}-\beta_{2a2}=0, \quad \beta_{2a1}+pq^{-1}\beta_{02}+\beta_{01}+pq^{-1}\beta_{2a2}=0.
\end{equation}
 It follows from (\ref{11_05_3}) that  $\beta_{01}=\beta_{02}=0$. Taking into account that $pq^{-1}\ne -1$ and the fact that $m$ is odd,   (\ref{22_04_3}) implies that $\beta_{a1}=\beta_{a2}=0$ for all $a\in \mathbb{Z}(m)$.
Equating the constant term on each side  of (\ref{22_04_2}) we receive
\begin{equation}\label{18_05_2}
\gamma_{2a1}+\gamma_{02}-\gamma_{01}-\gamma_{2a2}=0.
\end{equation}
Taking into account that in view of  (\ref{11_05_3}) $\gamma_{01}=\gamma_{02}=0$ and the
fact that $m$ is odd, (\ref{18_05_2}) implies that
  $\gamma_{a1}=\gamma_{a2}$ for all $a\in \mathbb{Z}(m)$.
Put $\gamma_{a1}=\gamma_{a2}=c_a$, $a\in \mathbb{Z}(m)$. Thus, we proved that the characteristic functions   $\hat M_j(s, a)$ have the following representation
\begin{equation}\label{22_04_4}
\hat M_j(s, a)=\exp\{-\sigma_j s^2+c_a\}, \quad (s, a)\in \mathbb{R}\times \mathbb{Z}(m),
\quad j=1, 2.
\end{equation}

It follows from   (\ref{19_05_1}) that the distributions
 $\nu_j$ are concentrated on the subgroup $\pi(\mathbb{R}\times \mathbb{Z}(mn))$. By Lemma \ref{lem3}, the distributions
 $\mu_j$ can be substituted by their  shifts $\mu'_j$ in such a manner that  $\nu_j=\mu'_j*\bar\mu'_j$ and the distributions
 $\mu'_j$ are also  concentrated on the subgroup $\pi(\mathbb{R}\times \mathbb{Z}(mn))$.
 Since $\pi$ is an isomorphism of the semigroups of distributions
     ${\rm M^1}(\mathbb{R}\times \mathbb{Z}(mn))$ and ${\rm M^1}(\pi(\mathbb{R}\times \mathbb{Z}(mn))),$ we have
  \begin{equation}\label{11_05_1}
M_j=\pi^{-1}(\nu_j)=\pi^{-1}(\mu'_j)*\pi^{-1}(\bar\mu'_j)=\pi^{-1}(\mu'_j)*
\overline{\pi^{-1}(\mu'_j)}, \quad j=1, 2.
\end{equation}
 Put $N_j=\pi^{-1}(\mu'_j)$, $j=1, 2$. We find from (\ref{11_05_1}) that
 \begin{equation}\label{11_05_4}
\hat M_j(s, a)=\hat N_j(s, a)\overline{\hat N_j(s, a)}, \quad (s, a)\in \mathbb{R}\times \mathbb{Z}(m), \quad j=1, 2.
\end{equation}
Putting in (\ref{11_05_4}) $a=0$  and taking into account (\ref{11_05_3}),
by the Cram\'er theorem on decomposition of the Gaussian distribution, we get
\begin{equation}\label{11_05_5}
\hat N_j(s, 0)=\exp\left\{-{\sigma_j\over 2} s^2+ib_{0j}s\right\}, \quad s\in \mathbb{R},
\quad j=1, 2,
\end{equation}
where $b_{0j}\in \mathbb{R}$. We will assume without loss of generality that  $b_{0j}=0$.

By Lemma  \ref{lem7},  $\hat N_j(s, a)$ are entire functions for any
$a\in \mathbb{Z}(m)$. It follows from  $(\ref{22_04_4})$ and  (\ref{11_05_4}) that they do
not vanish, and (\ref{8_16}) and  (\ref{11_05_5}) imply that the entire function
$\hat N_j(s, a)$ is of at most  order 2 and  type  ${\sigma_j\over 2}$. Applying the Hadamard theorem on the representation of an entire function of finite order,
(\ref{22_04_4}) and (\ref{11_05_4}), we obtain that
\begin{equation}\label{11_05_2}
\hat N_j(s, a)=\exp\left\{-{\sigma_j\over 2} s^2+ib_{aj}s+d_{aj}\right\}, \quad (s, a)\in \mathbb{R}\times \mathbb{Z}(m),  \quad j=1, 2,
\end{equation}
where $b_{aj}\in \mathbb{R}$. It follows from  (\ref{11_05_2}) that
$$
{\hat N_j(-iy+x, a)/\hat N_j(-iy, 0)}=
\exp\left\{-{\sigma_j\over 2} x^2+b_{aj}y+\sigma_jixy+b_{aj}ix+d_{aj}\right\}.
$$
By Lemma  \ref{lem7}, the function $
{\hat N_j(-iy+x, a)/\hat N_j(-iy, 0)}$
for any fixed $y$
is a characteristic function of variable $(x, a)\in\mathbb{R}\times \mathbb{Z}(m)$.  Thus, in particular, for any $y$ the inequality  $
|{\hat N_j(-iy+x, a)/\hat N_j(-iy, 0)}|\le 1$ holds true.
  It follows
from this that
$b_{aj}=0$, $j=1, 2,$ for any $a\in \mathbb{Z}(m)$, and (\ref{11_05_2})
implies that
\begin{equation}\label{11_05_6}
\hat N_j(s, a)=\exp\left\{-{\sigma_j\over 2} s^2+d_{aj}\right\},
\quad (s, a)\in \mathbb{R}\times \mathbb{Z}(m),  \quad j=1, 2.
\end{equation}

Let $P_j$ be the Gaussian distributions on $\mathbb{R}$ with the characteristic functions $\hat P_j(s)=\exp\{-{\sigma_j\over 2} s^2\}$, $s\in \mathbb{R}$. Denote by $Q_j$ the distributions on $\mathbb{Z}(m)$ with the characteristic functions $\hat Q_j(a)=\exp\{d_{aj}\}$, $a\in \mathbb{Z}(m)$. It follows from (\ref{11_05_6}) that $N_j=P_j*Q_j$, $j=1, 2$, and hence,  $\mu_j'=\pi(N_j)=\pi(P_j)*\pi(Q_j)$.
 It is obvious that $\pi(P_j)\in\Gamma(X)$. Verify that $\pi(\mathbb{Z}(m))=K$.
 We have
\begin{equation}\label{26_04_1}
A(Y, K)=(I+\tilde\alpha)(Y)=Y^{(m)}.
\end{equation}
Note that (\ref{23_04_1}) implies the equality
$$
Y^{(m)}=H\cup(my_0+H)\cup(2my_0+H)\cup\dots\cup((n-1)my_0+H),
$$
and hence, the equivalence
\begin{equation}\label{23_04_2}
y\in Y^{(m)}\Leftrightarrow\tau y\in \mathbb{R}\times \{0\}\times \mathbb{Z}(n).
\end{equation}
We have, $y\in A(Y, \pi(\mathbb{Z}(m)) \Leftrightarrow(g, \tau y)=1$ for all $g\in \mathbb{Z}(m)\Leftrightarrow\tau y\in \mathbb{R}\times \{0\}\times \mathbb{Z}(n).$ Taking
into account (\ref{26_04_1}) and (\ref{23_04_2}), this implies the equality
$A(Y, \pi(\mathbb{Z}(m))=A(Y, K),$ and hence, $\pi(\mathbb{Z}(m))=K$.
Thus, $\pi(Q_j)\in {\rm M}^1(K)$. So, we proved that $\mu'_j\in\Gamma(X)*{\rm M}^1(K)$, and
hence $\mu_j\in\Gamma(X)*{\rm M}^1(K)$, i.e. $\mu_j=\lambda_j*\rho_j$, where $\lambda_j\in \Gamma(X)$, $\rho_j\in {\rm M^1}(K)$, $j=1, 2.$

Let the characteristic functions    $\hat\lambda_j(y)$ be of the form
\begin{equation}\label{12_05_10}
\hat\lambda_j(y)=(x_j, y)\exp\{-\sigma_j y^2\}, \quad y\in Y\quad j=1, 2.
\end{equation}
Put $L=A(Y, K)$ and note that $\tilde\alpha(L)=L$. Since $\rho_j\in {\rm M^1}(K)$, we have $\hat\rho_j(y)=1$, $y\in L$, $j=1, 2$.
Consider the restriction of Heyde's functional equation (\ref{42}) for the characteristic
functions $\hat\mu_j(y)$  to $L$. Taking into account (\ref{12_05_10}), we find from the
obtained equation that
\begin{equation}\label{16_05_1}
\sigma_1+pq^{-1}\sigma_2=0,
\end{equation}
and hence,
\begin{equation}\label{12_05_12}
(x_1, u+v)(x_2, u+\tilde\alpha v)=(x_1, u-v)(x_2, u-\tilde\alpha v)\quad u, v\in L.
\end{equation}
It follows from (\ref{12_05_12}) that $2(x_1+\alpha x_2)\in K$. Since $f_2\in {\rm Aut}(X)$, we
have $x_1+\alpha x_2\in K$. Put $x_0=x_1+\alpha x_2$, $\gamma_1=\lambda_1*E_{-x_0}$, $\omega_1=\rho_1*E_{x_0}$, $\gamma_2=\lambda_2$, $\omega_2=\rho_2$. It is obvious that $\mu_j=\gamma_j*\omega_j$ and $\omega_j\in{\rm M}^1(K)$,  $j=1, 2$.
It is easy to see that the characteristic functions   $\hat\gamma_j(y)$ satisfy Heyde's
functional equation  (\ref{42}). Hence, the characteristic functions $\hat\omega_j(y)$
also satisfy  equation  (\ref{42}). We have
\begin{equation}\label{12_05_14}
\hat\omega_1(u+v)\hat\omega_2(u+\tilde\alpha v)=\hat\omega_1(u-v)\hat\omega_2(u-\tilde\alpha v),\quad u, v\in Y.
\end{equation}
We consider here $\omega_j$ as distributions on the group $X$. Denote by $M$ the character group
of the group $K$. Since $\omega_j\in{\rm M}^1(K)$,  $j=1, 2$, and $\alpha x=-x$ for any $x\in K$, we can consider
$\omega_j$ as distributions on the group $K$. Then equation (\ref{12_05_14})
takes the form
$$
\hat\omega_1(u+v)\hat\omega_2(u- v)=\hat\omega_1(u-v)\hat\omega_2(u+ v),\quad u, v\in M.
$$
Putting here $u=v=y$, we get $\hat\omega_1(2y)=\hat\omega_2(2y)$, $y\in M$. Since    $K\cong\mathbb{Z}(m)$, and $m$ is odd, this implies that $\hat\omega_1(y)=\hat\omega_2(y)$,
$y\in M$, and hence
$\omega_1=\omega_2=\omega$. Thus, $\mu_j=\gamma_j*\omega$, where $\gamma_j\in \Gamma(X)$, $\omega\in {\rm M^1}(K)$, $j=1, 2$.

Assume now that $\alpha>0$. Obviously, (\ref{16_05_1})  implies that   $\sigma_1=\sigma_2=0$,
i.e. $\gamma_j$ are degenerate distributions, and hence, the distributions $\mu_j$ are of the required form.
Thus, we proved Theorem 1 assuming that  $\alpha\ne\pm I$.

It still remains the cases when $\alpha=\pm I$. Note that since the group $X$ contains  no
elements of order $2$, we have $Y^{(2)}=Y$.
Let $\alpha=I$.
   Substituting in equation
(\ref{20_04_3}) $\alpha=I$, $u=v=y$, we receive $\hat\nu_1(2y)=\hat\nu_2(2y)=1$, $y\in Y$.
Hence, $\hat\nu_1(y)=\hat\nu_2(y)=1$, $y\in Y$, and this means that
$\nu_1=\nu_2=E_0$. Therefore,
$\mu_j$ are degenerate distributions. In this case Theorem 1 is also proved.
Let $\alpha=- I$. Then $K=X$. Substituting in equation
(\ref{42}) $\alpha=-I$, $u=v=y$, we obtain $\hat\mu_1(2y)=\hat\mu_2(2y)$, $y\in Y$.
Hence, $\hat\mu_1(y)=\hat\mu_2(y)$, $y\in Y$, and it means that
$\mu_1=\mu_2$.
Theorem 1 is completely proved.
$\Box$

{\bf Remark 1.} Consider an \text{\boldmath $a$}-adic
solenoid $X=\Sigma_\text{\boldmath $a$}$, and assume that $X$ contains no elements
of order $2$. Denote by  $Y$ the character group of the group $X$.
Let $\alpha = f_p f_q^{-1}$ be a topological automorphism of the group $X$.
 Put $K={\rm Ker}(I+\alpha)$. Theorem 1 can not be strengthened. Indeed, let $\gamma_j$ be Gaussian distributions on the group $X$ with the characteristic functions of the form
 $$\hat\gamma_j(y)=(x_j, y)\exp\{-\sigma_j y^2\},  \quad  y\in Y,
 $$ where  $x_1+\alpha x_2=0$ and $\sigma_1+pq^{-1}\sigma_2=0
$. Let $\omega\in {\rm M^1}(K)$. Put $\mu_j=\gamma_j*\omega$, $j=1, 2$. It easily follows from Lemma 1 that
if $\xi_1$ and $\xi_2$ are independent random variables with values in the group
       $X$  and distributions
  $\mu_1$ and $\mu_2$,  then the conditional distribution of the linear form
$L_2=\xi_1+\alpha\xi_2$
 given $L_1 = \xi_1 + \xi_2$ is symmetric. Thus, we can not narrow the class of distributions in Theorem 1 which is characterized by the symmetry of the conditional distribution of the linear form
$L_2$
 given $L_1$.

\bigskip

\centerline{\textbf{4.  Heyde's characterisation theorem for
\text{\boldmath $a$}-adic solenoids}}

\centerline{\textbf{containing an element of order $2$}}

\bigskip

 Let $X=\Sigma_{{\text{\boldmath $a$}}}$ and $G=X_{(2)}$.
 We discuss in this section the Heyde characterisation theorem for
\text{\boldmath $a$}-adic solenoids containing an element of order $2$, i.e. the case
when $G\cong \mathbb{Z}(2)$.
Let $\alpha$ be a topological automorphism of the group  $X$. Put $K={\rm Ker}(I+\alpha)$.
It is easy to see that $K\supset G$ for any $\alpha\in{\rm Aut}(X)$. Let $\xi_1$ and $\xi_2$
be independent random variables with values in the group $X$  and distributions
  $\mu_1$ and $\mu_2$ with nonvanishing characteristic functions. First we shall prove that generally speaking, Theorem 1 fails  if the group
  $X$ contains an element of order 2. Secondly, in the case when $K=G$ we give a complete description of distributions  which are characterized by the symmetry of the conditional distribution
 of the linear form    $L_2=\xi_1+\alpha\xi_2$ given $L_1 = \xi_1 + \xi_2$. It turns out that
 the corresponding class of distributions is wider than the class
$\Gamma(X)*{\rm M}^1(K)$.

Let $Y$ be the character group of the group $X$.  Note that a decomposition of the group $Y$
with respect to the subgroup
$Y^{(2)}$ consists of two cosets.
Let ${\mu\in {\rm M}^1(X)}$. It is easy to see that ${\mu\in \Gamma(X)*{\rm M}^1(G)}$ if and
only if  the characteristic function $\hat\mu(y)$ is represented in the form
$$
\hat\mu(y) = \begin{cases}(x, y)\exp\{-\sigma y^2\},
& y\in Y^{(2)},\\ (x, y)\kappa\exp\{-\sigma y^2\},
&   y\notin Y^{(2)},
\\
\end{cases}
$$
where $x\in X$,  $\sigma\ge 0$, $\kappa\in \mathbb{R}$, $|\kappa|\le 1$.
Introduce into consideration a class of distributions on the group
$X$ which is wider than the class   ${\Gamma(X)*{\rm M}^1(G)}$.

\bigskip

\noindent{\bf Definition  1}. {\it Let  $X=\Sigma_\text{\boldmath $a$}$,  $G=X_{(2)}.$ Assume
that $G\cong \mathbb{Z}(2)$. Let $Y$ be the character group of the group $X$. We say
that a distribution $\mu$ on the group $X$ belongs to the class $\Upsilon$  if its
characteristic function can be represented in the form
$$
\hat\mu(y) = \begin{cases}(x, y)\exp\{-\sigma y^2+i\beta y\}, & y\in Y^{(2)},\\ (x, y)\kappa\exp\{-\sigma' y^2+i\beta'y\}, &   y\notin Y^{(2)},
\\
\end{cases}
$$
for some $x\in X$,  $\sigma\ge 0$,  $\sigma'\ge 0$,  $\beta, \beta', \kappa\in \mathbb{R}$, $|\kappa|\le 1$.}

\bigskip

 Consider the group $\mathbb{R}\times \mathbb{Z}(2)$. Denote by $(t, k)$, $t\in \mathbb{R}$, $k\in \mathbb{Z}(2)$, its  elements. The character group of the group $\mathbb{R}\times \mathbb{Z}(2)$ is topologically isomorphic to the group $\mathbb{R}\times \mathbb{Z}(2)$.  Denote by
 $(s, l)$, $s\in \mathbb{R}$,
 $l\in \mathbb{Z}(2)$,
elements of the character group of the group $\mathbb{R}\times \mathbb{Z}(2)$.
Let $\mu\in \Gamma(\mathbb{R})*{\rm M}^1(\mathbb{Z}(2))$.
It is easy to see that the characteristic function of the distribution $\mu$ is of the form
$$
\hat\mu(s, l) = \begin{cases}\exp\{-\sigma s^2+i\beta s\},
&s\in \mathbb{R}, \ l=0,\\ \kappa\exp\{-\sigma s^2+i\beta s\},
&s\in \mathbb{R}, \   l=1,
\\
\end{cases}
$$
where $\sigma \ge 0$, $\beta, \kappa\in \mathbb{R}$, $|\kappa|\le 1$.
Introduce into consideration a class of distributions on the group
$X$ which is wider than $\Gamma(\mathbb{R})*{\rm M}^1(\mathbb{Z}(2))$.
 We need this class of distributions to study  distributions
 on \text{\boldmath $a$}-adic solenoids containing an element of order $2$
 which are characterized
by the symmetry of the conditional distribution of one linear form
 of   independent random variables given another.
\begin{lemma}\label{lem6}. {\it Consider the group  $\mathbb{R}\times \mathbb{Z}(2)$. Let $f(s, l)$ be a function on the character group of the group $\mathbb{R}\times \mathbb{Z}(2)$ of the form
\begin{equation}\label{y6}
f(s, l) = \begin{cases}\exp\{-\sigma s^2+i\beta s\}, &s\in \mathbb{R}, \ l=0,\\
\kappa\exp\{-\sigma' s^2+i\beta's\}, &s\in \mathbb{R}, \   l=1,
\\
\end{cases}
\end{equation}
where $\sigma\ge 0$,  $\sigma'\ge 0$,  $\beta, \beta', \kappa\in \mathbb{R}$. Then
$f(s, l)$ is the characteristic function of a signed measure $\mu$ on the group $\mathbb{R}\times \mathbb{Z}(2)$.
Moreover, $\mu$ is a measure if and only if either
$0<\sigma'<\sigma$ and $0<|\kappa|\le\sqrt{\sigma'\over \sigma}\exp\left\{-{(\beta-\beta')^2\over 4(\sigma-\sigma')}\right\}$
or $\sigma'=\sigma$, $\beta=\beta'$ and  $|\kappa|\le 1$. In the latter case $\mu\in \Gamma(\mathbb{R})*{\rm M}^1(\mathbb{Z}(2))$}.
\end{lemma}
\noindent{\bf Proof}. Obviously, we can prove the lemma assuming that $\kappa\ne0$. Multiplying if   necessary the function $f(s, l)$ by a suitable character
of the group $\mathbb{R}\times \mathbb{Z}(2)$, we can assume without loss of generality that $\kappa>0$.
Let $a\ge 0$. Denote by $\gamma_a$ the Gaussian distribution on the group
$\mathbb{R}$ with the characteristic function
\begin{equation}\label{17.01.16.1}
\hat\gamma_a(s)=\exp\{-a s^2\}, \quad s\in \mathbb{R}.
\end{equation}
Consider on the group $\mathbb{R}$ the measure
$\lambda_1={{1\over 2}(\gamma_{\sigma}*E_{\beta}
+\kappa\gamma_{\sigma'}*E_{\beta'})}$ and the signed measure  ${\lambda_2={1\over 2}(\gamma_{\sigma}*E_{\beta}
-\kappa\gamma_{\sigma'}*E_{\beta'})}$.
Define a signed measure $\mu$ on the group $\mathbb{R}\times \mathbb{Z}(2)$ in the
following way
$$\mu\{(E, 0)\}=\lambda_1\{E\}, \quad \mu\{(E, 1)\}=\lambda_2\{E\},
$$
where $E$ is a Borel subset of  $\mathbb{R}$. Taking into account that $\hat\lambda_1(s)+
\hat\lambda_2(s)=\hat\gamma_\sigma(s)\exp\{i\beta s\}$ and $\hat\lambda_1(s)-
\hat\lambda_2(s)=\kappa\hat\gamma_{\sigma'}(s)\exp\{i\beta's\}$, we have
$$
\hat\mu(s, l)=\int\limits_{\mathbb{R}\times
\mathbb{Z}(2)}e^{its}(k, l)d\mu(t, k)=\int\limits_{\mathbb{R}\times \{0\}}e^{its}d\mu(t, 0)+\int\limits_{\mathbb{R}\times \{1\}}e^{its}(1, l)d\mu(t, 1)=f(s, l).
$$
Thus, $f(s, l)$ is the characteristic function of the signed measure
 $\mu$, and the signed measure $\mu$ is a measure if and only if
   the signed measure $\lambda_2$ is a measure. It is obvious that
  if the signed measure $\lambda_2$ is a measure, then either  $\sigma > 0$ and  $\sigma'> 0$
  or $\sigma=\sigma'= 0$.
It is clear that if   $\sigma=\sigma'= 0$, then the signed measure
$\mu$ is a measure if and only if  $\beta=\beta'$ and  $\kappa\le 1$.
The statement of the lemma is proved in this case.

Let $\sigma > 0$ and  $\sigma'> 0$. Let $a>0$.   The density of the distribution
$\gamma_a$ is of the form
\begin{equation}\label{y5}
\rho_a(t)={1\over 2\sqrt{\pi a}}\exp\left\{-{t^2\over 4a}\right\},
\quad t\in \mathbb{R}.
\end{equation}
Taking into account (\ref{y5}),  the signed measure
$\lambda_2$ is a measure if and only if  the inequality
$$
{1\over 2\sqrt{\pi \sigma}}\exp\left\{-{(t-\beta)^2\over 4\sigma}\right\}-{\kappa\over 2\sqrt{\pi \sigma'}}\exp\left\{-{(t-\beta')^2\over 4\sigma'}\right\}\ge 0
$$
holds for all $t\in \mathbb{R}$.
This inequality is equivalent to the following
\begin{equation}\label{02_11_1}
\kappa\le\sqrt{\sigma'\over \sigma}\exp\left\{-{(t-\beta)^2\over 4\sigma}+{(t-\beta')^2\over 4\sigma'}\right\}, \quad t\in \mathbb{R}.
\end{equation}
 Let $\sigma\ne\sigma'$. Since $\kappa>0$, we have $\sigma'<\sigma$. The function in the right-hand side of inequality (\ref{02_11_1}) reaches its   minimum at $t={\sigma\beta'-\sigma'\beta\over \sigma-\sigma'}$, and this minimum is equal to  $\sqrt{\sigma'\over \sigma}\exp\left\{-{(\beta-\beta')^2\over 4(\sigma-\sigma')}\right\}$. Assume that $\sigma=\sigma'$. Then (\ref{02_11_1}) implies that $\beta=\beta'$ and  $ \kappa\le 1$. Thus,
 the signed measure
$\lambda_2$, and hence, the signed measure  $\mu$, is a measure if and only if
either  $\sigma'<\sigma$ and $d\le\sqrt{\sigma'\over \sigma}\exp\left\{-{(\beta-\beta')^2\over 4(\sigma-\sigma')}\right\}$
or $\sigma'=\sigma$, $\beta=\beta'$ and  $d\le 1$. It is also obvious that
in the latter case  $\mu\in \Gamma(\mathbb{R})*{\rm M}^1(\mathbb{Z}(2))$. Lemma 5 is proved.

We also note that if in  (\ref{y6}) $\sigma>0$ and $\kappa=
0$, then
$f(s, l)$ is the characteristic function of a
convolution of a non-degenerate Gaussian distribution on
 the group $\mathbb{R}$ and the Haar distribution on $\mathbb{Z}(2)$.
$\Box$

We   prove now that generally speaking, Theorem 1 fails  if the group
  $X=\Sigma_\text{\boldmath $a$}$ contains an element of order 2.

\bigskip

 \noindent{\bf Proposition 1.}  {\it  Let  $X=\Sigma_\text{\boldmath $a$}$,  $G=X_{(2)}.$ Assume that $G\cong \mathbb{Z}(2)$.
Let $\alpha$ be a topological automorphism of the group $X$, $\alpha<0$, $\alpha\ne -I$. Put $K={\rm Ker}(I+\alpha)$. Then
 there exist independent random variables $\xi_1$ and $\xi_2$ with values in the group $X$
 and  distributions $\mu_j\in \Upsilon$, $\mu_j\notin\Gamma(X)*{\rm M}^1(K)$, $j=1, 2$,  with
 nonvanishing characteristic functions such that the conditional distribution of the
 linear form $L_2 = \xi_1 + \alpha\xi_2$ given $L_1 = \xi_1 +
\xi_2$
is symmetric.}

\bigskip

\noindent{\bf Proof}.   Consider the group $\mathbb{R}\times \mathbb{Z}(2)$.
 Let $a\in {\rm Aut}(\mathbb{R}\times \mathbb{Z}(2))$. It is obvious that
 $a$ is of the form
 $a(t, k)=(c_a t, k)$, where $c_a\in \mathbb{R}$, $c_a\ne 0$.
 We will identify $a$ and $c_a$, i.e. we will write
$a(t, k)=(a t, k)$ and assume that
$a\in \mathbb{R}$, $a\ne 0$. The converse is also true. A nonzero real number
 corresponds to a topological automorphism of the group  $\mathbb{R}\times \mathbb{Z}(2)$.
 We remark that if $a\in {\rm Aut}(\mathbb{R}\times \mathbb{Z}(2))$, then
$\tilde a$ is of the form $\tilde a(s, l)=(a s, l)$, i.e. $\tilde a=a$.

Choose numbers  $\sigma_j$, $\sigma'_j$ and $\kappa$ in such a way that
 they satisfy the conditions:
 $0<\sigma_j'<\sigma_j$,  $0<\kappa\le \sqrt{\sigma_j'\over \sigma_j}$, $j=1, 2$.
Consider on the group $\mathbb{R}\times \mathbb{Z}(2)$ the functions
 \begin{equation}\label{y7}
f_j(s, l) = \begin{cases}\exp\{-\sigma_j s^2\}, &s\in \mathbb{R}, \ l=0,\\
\kappa\exp\{-\sigma_j' s^2\}, &s\in \mathbb{R}, \   l=1, \quad j=1, 2.
\\
\end{cases}
\end{equation}
By Lemma \ref{lem6}, there exist distributions $\lambda_j\in{\rm M}^1(\mathbb{R}\times \mathbb{Z}(2))$ such that
$\hat \lambda_j(s, l)=f_j(s, l)$. It is obvious that $\lambda_j\notin\Gamma(\mathbb{R})*{\rm M}^1(\mathbb{Z}(2))$, $j=1, 2$. By the condition $\alpha = f_p f_q^{-1}$ for some mutually prime  $p$ and $q$,
where $f_p, f_q \in {\rm
Aut}(X)$. It follows from $G\cong \mathbb{Z}(2)$, that $p$ and $q$ are odd,
and hence $f_p$ and $f_q$ are topological automorphisms of the group  $\mathbb{R}\times \mathbb{Z}(2)$.
Since $\alpha<0$, we can assume that the conditions
\begin{equation}\label{y8}
\sigma_1+pq^{-1}\sigma_2=0,
\end{equation}
\begin{equation}\label{y9}
\sigma'_1+pq^{-1}\sigma'_2=0
\end{equation}
hold.
Verify that the characteristic functions $\hat \lambda_j(s, l)$ satisfy Heyde's functional
equation  $(\ref{42})$.
Put $u=(s_1, l_1)$, $v=(s_2, l_2)$. Let either $u, v\in \mathbb{R}\times \{0\}$ or $u, v\in \mathbb{R}\times \{1\}$. Then $u\pm v, u\pm \alpha v\in \mathbb{R}\times \{0\}$,
and taking into account (\ref{y7}),  equality $(\ref{42})$ for these $u, v$ is satisfied if
the equality
$$
\sigma_1(s_1+s_2)^2+\sigma_2(s_1+pq^{-1} s_2)^2=\sigma_1(s_1-s_2)^2+\sigma_2(s_1-pq^{-1} s_2)^2, \quad s_1, s_2\in \mathbb{R},
$$
holds true. But this  equality follows from $(\ref{y8})$. Let
either $u\in \mathbb{R}\times \{0\}$, $v\in \mathbb{R}\times \{1\}$ or $u\in \mathbb{R}\times \{1\}$, $v\in \mathbb{R}\times \{0\}$. Then $u\pm v, u\pm \alpha v\in \mathbb{R}\times \{1\}$,
and taking into account (\ref{y7}),  equality $(\ref{42})$ for these $u, v$
is satisfied if  the equality
$$
\sigma'_1(s_1+s_2)^2+\sigma'_2(s_1+pq^{-1} s_2)^2=
\sigma'_1(s_1-s_2)^2+\sigma'_2(s_1-pq^{-1} s_2)^2, \quad s_1, s_2\in \mathbb{R},
$$
is valid. But this equality follows from $(\ref{y9})$. It is obvious that we exhausted all possibilities for    $u$ and $v$. Thus, we verified that the characteristic functions
$\hat\lambda_j(l, s)$ satisfy Heyde's functional equation $(\ref{42})$.

Let  $\zeta_1$ and  $\zeta_2$ be independent random variables with values
in the group $\mathbb{R}\times \mathbb{Z}(2)$  and distributions $\lambda_1$ and $\lambda_2$.
Since the characteristic functions $\hat\lambda_j(l, s)$ satisfy Heyde's functional  equation (\ref{42}), by Lemma \ref{lem1},
the conditional distribution of the linear form
$T_2=\zeta_1+\alpha\zeta_2$ given $T_1=\zeta_1+\zeta_2$ is symmetric.

Denote by  $Y$ the character group of the group  $X$. Let $y_0\notin Y^{(2)}$ and $Y=Y^{(2)}\cup(y_0+Y^{(2)})$
be a decomposition of the group $Y$ with respect to the subgroup $Y^{(2)}$.
Define the mapping
 $\tau: Y\mapsto\mathbb{R}\times \mathbb{Z}(2)$ by the formula
  \begin{equation}\label{5_16}
\tau(y)= \begin{cases}(y, 0), & y\in Y^{(2)},\\ (y, 1), & y\in  y_0+Y^{(2)}.
\\
\end{cases}
\end{equation}
It is obvious that  $\tau$ is a homomorphism.  Put $\pi=\tilde\tau$,
$\pi: \mathbb{R}\times \mathbb{Z}(2)\mapsto X$.
Since the subgroup $\tau(Y)$ is dense in $\mathbb{R}\times \mathbb{Z}(2)$,
$\pi$ is a monomorphism.
   It is easy to see that
\begin{equation}\label{07.06.15.1}
\alpha\pi(t, k)=\pi f_p f_q^{-1}(t, k), \quad (t, k)\in \mathbb{R}\times \mathbb{Z}(2).
\end{equation}
  Put $\xi_j=\pi(\zeta_j)$,
$j=1, 2$. Then $\xi_j$ are independent random variables with values in
the group
$X$ and distributions $\mu_j=\pi(\lambda_j)$.
We have $\hat\mu_j(y)=\hat\lambda_j(\tau (y))$,  and (\ref{y7}) and (\ref{5_16})   imply that
\begin{equation}\label{07.01.16.1}
\hat\mu_j(y) = \begin{cases}\exp\{-\sigma_j y^2\}, & y\in Y^{(2)},\\ \kappa\exp\{-\sigma_j' y^2\}, &   y\in  y_0+Y^{(2)}, \quad j=1, 2.
\\
\end{cases}
\end{equation}
Obviously,  the conditional distribution of the linear form
 $L_2=\pi(T_2)$ given $L_1=\pi(T_1)$  is symmetric.
We have $L_1=\xi_1+\xi_2$, and it follows from (\ref{07.06.15.1}) that
$L_2=\xi_1+\alpha\xi_2$.  Definition 1 and  (\ref{07.01.16.1}) imply that $\mu_j\in \Upsilon$. Since $\sigma_j\ne \sigma'_j$ and $\alpha\ne -I$, we have  $\mu_j\notin \Gamma(X)*{\rm M}^1(K)$, $j=1, 2$. Proposition 1 is proved. $\Box$
\bigskip

\noindent{\bf Remark 2.} We consider here the  case  when $\alpha=-I$ and complete
Proposition 1 by  the following  statement.

{\it  Let  $X=\Sigma_\text{\boldmath $a$}$, $G=X_{(2)}$ and  $G\cong \mathbb{Z}(2).$
Let $\xi_1$ and $\xi_2$ be independent random variables
with values in the group
       $X$  and distributions
  $\mu_1$ and $\mu_2$ with nonvanishing characteristic functions.
  The conditional distribution of the linear form
$L_2=\xi_1-\xi_2$
 given $L_1 = \xi_1 + \xi_2$ is symmetric if and only if   either $\mu_1=\mu_2*\lambda$ or   $\mu_2=\mu_1*\lambda$, where  $\lambda\in {\rm M}^1(G)$.}

  Let $Y$ be the character group of the group $X$,  $y_0\notin Y^{(2)}$ and let $Y=Y^{(2)}\cup(y_0+Y^{(2)})$
be a decomposition of the group $Y$ with respect to the subgroup $Y^{(2)}$.
Assume that the conditional distribution of the linear form
$L_2$
 given $L_1$ is symmetric.    By Lemma \ref{lem1}, the characteristic functions $\hat\mu_j(y)$  satisfy Heyde's functional equation   $(\ref{42})$ which takes the form
\begin{equation}\label{18_08_1}
\hat\mu_1(u+v )\hat\mu_2(u-v)=
\hat\mu_1(u-v)\hat\mu_2(u+v), \quad u, v \in Y.
\end{equation}
 Substituting $u=v=y$ in equation $(\ref{18_08_1})$, we get $\hat\mu_1(2y)=\hat\mu_2(2y)$, $y\in Y$. Hence,
\begin{equation}\label{18_08_2}
\hat\mu_1(y)=\hat\mu_2(y), \quad y\in Y^{(2)}.
\end{equation}
Substituting $u=y$, $v=y+y_0$ in equation $(\ref{18_08_1})$, we obtain
\begin{equation}\label{4_11_1}
\hat\mu_1(2y+y_0)\hat\mu_2(y_0)=
\hat\mu_1(y_0)\hat\mu_2(2y+y_0), \quad y\in Y.
\end{equation}
Assume that ${|\hat\mu_1(y_0)|\le
|\hat\mu_2(y_0)|}$. We find from (\ref{4_11_1}) that
\begin{equation}\label{18_08_4}
{\hat\mu_1(2y+y_0)\over\hat\mu_2(2y+y_0)}={\hat\mu_1(y_0)\over
\hat\mu_2(y_0)}, \quad y\in Y.
\end{equation}
Put $d={\hat\mu_1(y_0)/
\hat\mu_2(y_0)}$. It follows from $\hat\mu_j(-y)=\overline{\hat\mu_j(y)}$, $y\in Y$, that $d$ is a real number, and (\ref{18_08_4}) implies that
\begin{equation}\label{18_08_5}
\hat\mu_1(y)=d\hat\mu_2(y), \quad y\in y_0+Y^{(2)}.
\end{equation}
Denote by $\lambda$ a distribution  on $G$ with the characteristic function
\begin{equation}\label{18_08_6}
\hat\lambda(y) = \begin{cases}1, & y\in Y^{(2)},\\ d, &   y\in y_0+Y^{(2)}.
\\
\end{cases}
\end{equation}
It follows from (\ref{18_08_2}), (\ref{18_08_5}) and
 (\ref{18_08_6}) that $\hat\mu_1(y)=\hat\mu_2(y)\hat\lambda(y),$ $y\in Y$, and hence $\mu_1=\mu_2*\lambda$.

 The converse statement follows directly from Lemma  \ref{lem1}.

\bigskip

In the case when $K=G$ and $G\cong \mathbb{Z}(2)$ we give a complete description of distributions
which are characterized by the symmetry of the conditional distribution
 of the linear form    $L_2$ given $L_1$.

\bigskip

\noindent{\bf Theorem 2.}  {\it  Let  $X=\Sigma_\text{\boldmath $a$}$,  $G=X_{(2)}.$
Let $\alpha$ be a topological automorphism of the group $X$. Put $K={\rm Ker}(I+\alpha)$. Assume that $K=G$ and $G\cong \mathbb{Z}(2)$. Let $\xi_1$ and $\xi_2$ be independent random variables
with values in the group
       $X$  and distributions
  $\mu_1$ and $\mu_2$ with nonvanishing characteristic functions.
  Assume that the conditional distribution of the linear form
$L_2=\xi_1+\alpha\xi_2$
 given $L_1 = \xi_1 + \xi_2$ is symmetric.
 Then $\mu_j\in \Upsilon $, $j=1, 2$. Moreover, if $\alpha>0$, then some sifts
of the distributions   $\mu_j$
are supported in $G$.}

\bigskip

 To prove Theorem 2 we need such lemmas.

\begin{lemma}\label{lem4} {\rm (\cite{My1})}. {\it Let  $X$ be a locally compact Abelian group, $\alpha$
be a topological automorphism of $X$. Let  $\xi_1$ and  $\xi_2$ be
independent random variables with values in
the group $X$.
   If the conditional distribution of the linear form    $L_2 = \xi_1 +
\alpha\xi_2$  given $L_1 = \xi_1 +
\xi_2$  is symmetric,   then the linear forms
$L_1'=(I+\alpha)\xi_1+2\alpha\xi_2$ and
$   L_2'=2\xi_1+(I+\alpha)\xi_2$ are independent.}
\end{lemma}

It is well-known that the Gaussian distribution on the real line
is characterized by the independence of two linear forms
of $n$ independent random variables. We need a group analogue of this theorem for
 $n=2$.
\begin{lemma}\label{lem5}{\rm (\cite{Fe16}, see also \cite[Theorem 10.3]{Fe5})}. \textit{Let  $X$ be a locally compact Abelian group containing
  no subgroup topologically isomorphic to the circle group
    $\mathbb{T}$. Let $\alpha_j$, $\beta_j$ be topological automorphisms
    of the group  $X$. Assume that $\xi_1$ and  $\xi_2$ are independent random variables with values in the group
$X$   and distributions   $\mu_1$ and $\mu_2$ with nonvanishing characteristic functions.
Then the independence of the linear forms $L_1 = \alpha_1\xi_1 + \alpha_2\xi_2$
and $L_2 = \beta_1\xi_1 +
\beta_2\xi_2$ implies that
   $\mu_j\in
\Gamma(X)$, $j=1, 2$.}
\end{lemma}

We will prove now that the characteristic functions of  divisors of a distribution on the group $  \mathbb{R}\times\mathbb{Z}(2)$ with the nonvanishing characteristic function  of the form (\ref{y6})  have the similar representation.
\begin{lemma}\label{lem8}. {\it Let   $N\in {\rm M}^1(\mathbb{R}\times
 \mathbb{Z}(2))$ and the characteristic function
  $\hat N(s, l)=f(s, l)$ is represented in the form
  $(\ref{y6})$, where  either
$0<\sigma'<\sigma$ and $0<|\kappa|\le\sqrt{\sigma'\over \sigma}\exp\left\{-{(\beta-\beta')^2\over 4(\sigma-\sigma')}\right\}$
or $\sigma'=\sigma$, $\beta=\beta'$ and  $0<|\kappa|\le 1$.
 Let
$N=P_1*P_2,$ where $P_j\in {\rm M}^1(\mathbb{R}\times \mathbb{Z}(2))$. Then each of the characteristic functions $\hat P_j(s, l)$ is of the form
$$
\hat P_j(s, l) = \begin{cases}\exp\{-\sigma_j s^2+i\beta_j s\}, &s\in \mathbb{R}, \ l=0,\\ \kappa_j\exp\{-\sigma_j' s^2+i\beta_j's\},
&s\in \mathbb{R}, \   l=1,
\\
\end{cases}
$$
 where either $0<\sigma_j'<\sigma_j$ and $0<|\kappa_j|\le\sqrt{\sigma_j'\over \sigma_j}\exp\left\{-{(\beta_j-\beta_j')^2\over 4(\sigma_j-\sigma_j')}\right\}$,
or $\sigma_j'=\sigma_j$, $\beta_j=\beta_j'$ and  $0<|\kappa_j|\le 1$.}
\end{lemma}
{\bf Proof.} We have
\begin{equation}\label{9_16}
\hat N(s, l)=\hat P_1(s, l)\hat P_2(s, l), \quad s\in \mathbb{R}, \ l\in \mathbb{Z}(2).
\end{equation}
Putting in (\ref{9_16}) $l=0$, we obtain  from  $(\ref{y6})$ by the Cram\'er theorem
on decomposition of the Gaussian distribution on the real line,
\begin{equation}\label{10_16}
\hat P_j(s, 0) = \exp\{-\sigma_j s^2+i\beta_js\}, \quad s\in \mathbb{R},
\quad j=1, 2,
\end{equation}
where $\sigma_j\ge 0$,     $\beta_j \in \mathbb{R}$.
By Lemma \ref{lem7},  $\hat P_j(s, 1)$ are entire functions, and $(\ref{y6})$ and  (\ref{9_16}) imply that they do not vanish. It follows from (\ref{8_16}) and (\ref{10_16}) that  the entire function
$\hat P_j(s, 1)$ is of at most order 2 and  type $\sigma_j$.
Taking into account that $\hat P_j(-s, 1)=\overline{\hat P_j(s, 1)}$, by the Hadamard theorem on the representation of an entire function
 of finite order and Lemma \ref{lem6},  we obtain the desired statement. Lemma \ref{lem8} is proved. $\Box$

\bigskip

\noindent{\bf Proof of Theorem 2.} Denote by $Y$ the character group of the group $X$ and consider  $Y=H_{\text{\boldmath $a$}}$  as a subset of  $\mathbb{R}$. Let $\alpha = f_p f_q^{-1}$ for some mutually prime $p$ and $q$,
where $f_p, f_q \in {\rm
Aut}(X)$. It follows from $G\cong \mathbb{Z}(2)$  that $p$ and $q$ are odd. Consider the distributions   $\nu_j = \mu_j* \bar \mu_j$ and the functions $\varphi_j(y) =  \log \hat\nu_j(y)$, $j=1, 2$. Reasoning as in the proof of Theorem 1, we see that the functions $\varphi_j(y)$ satisfy equation  (\ref{17}), where $H$ is of the form (\ref{20_08_1}).
Since $G\cong \mathbb{Z}(2)$ and
$K=G$, it is easy to see that
\begin{equation}\label{1_16}
I+\alpha=f_2\beta,
\end{equation}
where $\beta \in {\rm Aut}(X)$.
Since $p$ and  $q$ are odd,  (\ref{20_08_1}) and (\ref{1_16}) imply that
$H=(I-\tilde\alpha)(Y)=Y^{(n)}$, where $n$ is even, and we can assume that $n$ is  the minimum possible.
Consider the factor-group $Y/H$. It is obvious that $Y/H\cong \mathbb{Z}(n)$.  Take an element
 $y_0\in Y$ in such a way that the coset $y_0+H$  is a generator of the factor-group
  $Y/H$.
Let
$$
Y=H\cup(y_0+H)\cup(2y_0+H)\cup\cdots\cup((n-1)y_0+H)
$$
be a decomposition of the group $Y$ with respect to the subgroup $H$.
 It follows from (\ref{17}) that for any coset $ly_0+H$ there exist polynomials $P_{lj}(y)$
 on $\mathbb{R}$ with real coefficients of the form (\ref{13_05_1}) such that
 \begin{equation}\label{13_05}
P_{lj}(y)=\varphi_j(y), \quad y\in ly_0+H, \quad l=0, 1,\dots, n-1, \quad j=1, 2.
\end{equation}

  By Lemma \ref{lem4}, the linear forms
$  L_1'=(I+\alpha)\xi_1+2\alpha\xi_2$ and
$ L_2'=2\xi_1+(I+\alpha)\xi_2$ are independent. Consider the new independent random variables $\eta_j=2\xi_j$, $j=1, 2$.
   Taking into account (\ref{1_16}),  the linear forms $L_1'$ and $L_2'$ can be written as follows
$L_1'=\beta\eta_1+\alpha\eta_2$ and  $L_2'=\eta_1+\beta\eta_2$.
Since $\alpha, \beta\in {\rm Aut}(X)$ and the group  $X$ contains no subgroup topologically isomorphic to the circle group $\mathbb{T}$, by Lemma \ref{lem5},  the random variables $\eta_j$ have Gaussian distributions. This easily implies that
\begin{equation}\label{2_16}
\varphi_j(y)=-\sigma_jy^2, \quad y\in Y^{(2)}, \quad j=1, 2,
\end{equation}
where $\sigma_j\ge 0$.

  It is obvious that
$
Y=Y^{(2)}\cup(y_0+Y^{(2)})
$
is a decomposition of the group $Y$ with respect to the subgroup $Y^{(2)}$.
It follows from  (\ref{2_16}) and the inequality
$$
|\hat\nu_j(u)-\hat\nu_j(v)|^2\le 2(1-{\rm Re}\
\hat\nu_j(u-v)), \quad u, v\in Y, \quad j=1, 2,
$$
which holds for any characteristic function on the group  $Y$,
 that each of the functions  $\hat\nu_j(y)$ is uniformly continuous on
the coset $y_0+Y^{(2)}$ in the topology induced on
$y_0+Y^{(2)}$ by the
topology of the group  $\mathbb{R}$. Since
$$
y_0+Y^{(2)}=(y_0+H)\cup(3y_0+H)\cup\dots\cup((y_0-1)c+H)
$$
and each coset
$ly_0+H$ is a dense subset of $\mathbb{R}$, taking into account (\ref{13_05}) this implies that there exist polynomials
$A_jy^2+B_jy+C_j$,  $j=1, 2,$  such that
$$
P_{lj}(y)=A_jy^2+B_jy+C_j, \quad y\in ly_0+H, \quad l=1, 3, \dots, n-1, \quad j=1, 2.
$$
and hence, the representations
\begin{equation}\label{3_16}
\varphi_j(y)=A_jy^2+B_jy+C_j, \quad y\in y_0+Y^{(2)}, \quad  j=1, 2,
\end{equation}
hold. Since
  $\varphi_j(-y)=\varphi_j(y)$, $y\in Y$, we have $B_j=0,$ $j=1, 2$.
Put $A_j=-\sigma_j'$, $\kappa_j=e^{C_j}$.
As a result, (\ref{2_16})  and (\ref{3_16}) imply the representations
\begin{equation}\label{4_16}
\hat\nu_j(y) = \begin{cases}\exp\{-\sigma_j y^2\}, & y\in Y^{(2)},\\ \kappa_j\exp\{-\sigma_j' y^2\}, &   y\in y_0+Y^{(2)},
\\
\end{cases}
\end{equation}
where $0<\kappa_j\le 1$, $\sigma_j\ge 0$, $\sigma_j'\ge 0$, $j=1, 2$.

By Lemma 5, there exist signed measures $M_j$ on the group $\mathbb{R}\times \mathbb{Z}(2)$   with the characteristic functions of the form
\begin{equation}\label{6_16}
\hat M_j(s, l)=\begin{cases}\exp\{-\sigma_j s^2\}, &s\in \mathbb{R}, \ l=0,\\  \kappa_j\exp\{-\sigma_j' s^2\}, &s\in \mathbb{R}, \   l=1, \quad j=1, 2.
\\
\end{cases}
\end{equation}

Let a homomorphism $\tau: Y\mapsto\mathbb{R}\times \mathbb{Z}(2)$ is defined by formula (\ref{5_16}),  and $\pi=\tilde\tau$,
$\pi: \mathbb{R}\times \mathbb{Z}(2)\mapsto X$.
  Taking into account (\ref{5_16}), it follows from (\ref{4_16}) and (\ref{6_16}) that $\nu_j=\pi(M_j)$, $j=1, 2$.
Since $\pi$ is a continuous monomorphism, the signed measures $M_j$  are measures.
Obviously, the distributions $\nu_j$ are concentrated on the Borel subgroup
$F=\pi(\mathbb{R}\times \mathbb{Z}(2))$. By Lemma \ref{lem3}, the distributions $\mu_j$ can be substituted by their  shifts $\mu_j'$ in such a way that the distributions $\mu_j'$ are also concentrated on $F$ and $\nu_j = \mu_j'* \bar \mu_j'$, $j=1, 2$. Put $N_j=\pi^{-1}(\mu_j')\in {\rm M}^1(\mathbb{R}\times \mathbb{Z}(2))$.
 It is obvious that $M_j=N_j*\overline N_j$.
By Lemma  \ref{lem8},  the characteristic functions $\hat N_j(s, l)$  are of the form
\begin{equation}\label{14_16}
\hat N_j(s, l)=\begin{cases}\exp\{-{\sigma_j \over 2}s^2+i\beta_js \}, &s\in \mathbb{R}, \ l=0,\\ \sqrt{\kappa_j}\exp\{-{\sigma'_j \over 2}s^2+i\beta'_js\}, &s\in \mathbb{R}, \   l=1,
\\
\end{cases}
\end{equation}
where $\beta_j, \beta_j'\in \mathbb{R}$,  $j=1, 2$.
Since $\mu_j'=\pi(N_j)$, we have $\hat\mu_j'(y)=\hat N_j(\tau(y))$ and it follows from   $(\ref{5_16})$ and (\ref{14_16})  that the characteristic functions $\hat\mu_j'(y)$  are of the form
$$
\hat \mu'_j(y)=\begin{cases}\exp\{-{\sigma_j \over 2}y^2+i\beta_jy \}, & y\in Y^{(2)},\\ \sqrt{\kappa_j}\exp\{-{\sigma'_j \over 2}y^2+i\beta'_jy\}, &   y\in y_0+Y^{(2)}, \quad j=1, 2.
\\
\end{cases}
$$
Hence, $\mu_j\in \Upsilon$, $j=1, 2$.

Let $\alpha>0$.  Assuming that $u, v\in Y^{(2)}$, substitute (\ref{2_16}) in (\ref{12}). We obtain $\sigma_1+pq^{-1}\sigma_2=0$, and this implies that   $\sigma_1=\sigma_2=0$. Hence, $\hat\nu_j(y)=1$, $y\in Y^{(2)}$.
It follows from  Lemma \ref{lem2} that $\sigma(\nu_j)\subset A(X, Y^{(2)})=X_{(2)}=G$, $j=1, 2$. It means by Lemma \ref{lem3}, that some shifts of the distributions $\mu_j$
are supported in $G.$
Theorem 2 is completely proved. $\Box$

\end{document}